\documentstyle{amsppt}
\magnification=1200

\topmatter

\title On Chevalley restriction theorem\endtitle

\author E. A. Tevelev\endauthor

\endtopmatter

\def\g{\frak g}
\def\a{\frak a}
\def\h{\frak h}
\def\id{\text{id}}
\def\C{\Bbb C}
\def\R{\Bbb R}
\def\gR{\g^\R}
\def\k{\frak k}
\def\b{\frak b}
\def\kR{\k^\R}
\def\aR{\a^\R}
\def\p{\frak p}
\def\pR{\p^\R}
\def\isom{\simeq}
\def\Ker{\mathop{\text{\rm Ker}}}
\def\Lie{\mathop{\text{\rm Lie}}}
\def\Im{\mathop{\text{\rm Im}}}
\def\Supp{\mathop{\text{\rm Supp}}}

\document

\head \S0. Introduction\endhead
Let $\g$ be a complex semisimple Lie algebra with adjoint group $G$.
Suppose that $\sigma\ne\id$ is an involutive automorphism of $\g$.
Then $\sigma$ induces uniquely an involution of $G$ also denoted by $\sigma$,
let $K=G^\sigma$ be a subgroup of $\sigma$-fixed points.
Consider a direct decomposition $\g=\k\oplus\p$
of $\g$ into eigenspaces for $\sigma$, so
$$\k=\{x\in\g\,|\,\sigma(x)=x\},\qquad \p=\{x\in\g\,|\,\sigma(x)=-x\}.$$
Then, clearly, $\k$ is a Lie algebra of $K$ and $\p$ is a $K$-module,
this representation is certainly nothing else but an isotropy
representation of a symmetric space $G/K$
(see e.g.~\cite{KR}).
Denote by $\a\subset\p$ any maximal abelian ad-diagonalizable subalgebra.
Consider the ``baby Weyl group'' $W_0=N_K(\a)/Z_K(\a)$.
It is well-known that $W_0$ is a finite group
generated by reflections as a linear group operating on $\a$
(with respect to some real form of $\a$).
Let $\psi:\,\C[\p]^K\to\C[\a]^{W_0}$ be a restriction map of algebras of invariants.
Then the famous Chevalley restriction theorem
states that $\psi$ is an isomorphism (see eg.~\cite{He1},
in \cite{Vi} a more general result is obtained
in the context of so-called $\theta$-groups
attached to any periodic automorphism of $\g$).

Now consider the following special case.
The adjoint representation $G:\g$ could be identified
with an isotropy representation of a symmetric space
induced by an involution $\sigma$ on $G\times G$ given by $\sigma(x,y)=(y,x)$.
Then the cited Chevalley restriction theorem
is equivalent to the ``usual'' Chevalley restriction theorem
$\C[\g]^G\isom\C[\h]^W$, where $\h\subset \g$ is a Cartan subalgebra
and $W$ is the usual Weyl group.
In the paper~\cite{J} A.~Joseph obtained the remarkable
``multivariable'' analogue of this theorem.
Namely, he proved that the restriction map
$$\C[\g\times\g]^G\to\C[\h\times\h]^W$$
is surjective.
This theorem has important applications to
the representation theory and the geometry of commuting varieties
(cf.~\cite{Ri}), because it is essentially equivalent to
the fact that the ring of $G$-invariant functions on a commuting
variety is integrally closed (see \cite{Hu, J}).
The aim of this paper is to extend Joseph's arguments in order to
prove the following generalization:

\proclaim{Theorem}
The restriction map $\psi:\,\C[\p\times\p]^K\to\C[\a\times\a]^{W_0}$
is surjective.
\endproclaim

Actually, this (and Joseph's) theorem
also holds for any number of summands (generalization of the proof is immediate),
but we will prove the Theorem in this form in order to make the notation easier.

Certainly, one could try to deduce this Theorem directly from Joseph's
result arguing as follows:
take $f\in\C[\a\times\a]^{W_0}$, lift it to an element
$\C[\h\times\h]^W$, then by Joseph's theorem this element
could be lifted to an element of $\C[\g\times\g]^G$ and we could resrict it
to an element of $\C[\p\times\p]^{K}$.
Unfortunately, the first step of this construction fails
even in the case of one summand (see examples in \cite{He2}).
Moreover, I don't know how to prove the natural
analogue of the Theorem in the context of arbitrary $\theta$-groups
(see \cite{Vi}) and this seems to be an interesting problem.
I should mention also that in some particular cases
these results could be proved using
``generalized polarizations"  (see eg.~\cite{Hu}), but these results
don't cover our Theorem (as well as Joseph's).
In \cite{Pa1} Theorem was proved in the case
of an involution of maximal rank.

The paper is organized as follows.
In \S1 we will look for a proof of Chevalley
restriction theorem in the case of one summand such that this proof could be
generalized to the case of two summands.
Another reason for doing that separately
is that all the ideas distinct from the ones used by Joseph
will already appear in this case.
In \S2 we will prove the Theorem
in full generality showing that our main tool (class~$1$ representations)
is compatible with Joseph's idea to use the PRV Conjecture
(or, better said, the Kumar-Mathieu Theorem) to ``separate layers".

I would like to thank D.~Saltman for useful discussions
and his warm hospitality during my stay in the University
of Texas in Austin. The research was supported by CRDF grant RM1-206,
and by INTAS grant INTAS-OPEN-97-1570.

\head\S1. Surjectivity of $\psi:\,\C[\p]^K\to\C[\a]^{W_0}$\endhead
Recall, that there is 1-1 correspondence between involutions of $\g$
and real forms of $\g$ (see \cite{W, \S1.1}).
Namely, suppose that $\gR$ is a real form of $\g$
(that is, $\g=\gR+i\gR$), let $\gR=\kR\oplus\pR$ be a Cartan decomposition,
let $\g=\k\oplus\p$ be its complexification.
Then the map $\sigma$ defined by $\sigma(x)=x$ for any $x\in\k$ and
$\sigma(x)=-x$ for any $x\in\p$ is an involution
and any involution can be obtained this way.
Now, a maximal ad-diagonalizable abelian subalgebra $\a$ in $\p$
can be chosen to be a complexification of some $\aR\subset\pR$.
From now on by complex conjugation we will understand the complex conjugation
with respect to a given real form $\gR$.

Let $B$ be a Killing form on $\g$.
Then the restriction of $B$ on $\kR$ is negative-definite
and on $\pR$ is positive definite.
In particular we can identify $\g\isom\g^*$ and $\a\isom\a^*$ by means
of this scalar product. Then the restriction map
$\psi:\,\C[\p]^K\to\C[\a]^{W_0}$
can be viewed as a homomorphism of symmetric algebras
$\psi:\,(S\p)^K\to (S\a)^{W_0}$.
Now, it is well-known that $K=G^\sigma$ is a reductive group
(see eg.~\cite{Vu}), so we have a $K$-module decomposition
$S\p=(S\p)^K\oplus(S\p)_K$, where
$(S\p)_K$ is a sum of all non-trivial $K$-submodules.
Consider the Reynolds operator
($K$-equivariant projection on the first factor of this decomposition).
Let $\theta:\,(S\a)^{W_0}\to(S\p)^K$
be its restriction to $(S\a)^{W_0}$.
It is clear that we have a commutative diagram
$$\matrix
S\a&\mathop{\longrightarrow}\limits^{\theta_1}& (S\p)^K\cr
{}_{\theta_2}\searrow\!\!\!\!\!\!\!{}&&{}\!\!\!\!\!\!\!\!\!\!\!\!\!\!\nearrow_{\theta}\cr
&(S\a)^{W_0}&\cr\endmatrix$$
where $\theta_1$ is a restriction of the Reynolds operator to $S\a$
and $\theta_2$ is a $W_0$-equivariant projection.
Now we are ready to prove the following lemma:

\proclaim{Lemma 1.1}
If $\theta$ is injective then $\psi$ is surjective.
\endproclaim

\demo{Proof}
It is clear that both maps $\psi$ and $\theta$ are homogeneous
(with respect to natural grading of $S\a$ and $S\p$).
So the composition $\psi\theta$ is a homogeneous map of degree~$0$
from $(S\a)^{W_0}$ to $(S\a)^{W_0}$.
Now it is clear that if $\theta$ is injective and $\Ker\psi\cap\Im\theta=0$
then $\psi\theta$ is injective, hence surjective, so $\psi$ is surjective as well.
Certainly the second assertion follows from the fact that $\psi$
is injective (and this is an easy part of Chevalley restriction theorem),
but we are looking for a proof that could be applicable in the case
of two summands.

First we could extend $B$ to the scalar product in $S\p$ in the usual way.
Since the restriction of $B$ on $\pR$ is positive definite, it is clear that
the restriction of $B$ on $S\pR$ is also positive definite (more precisely,
the restriction of $B$ on any finite-dimensional subspace of $S\pR$
is positive definite). In particular, if $V\subset S\p$ is any subspace
stable under complex conjugation, then $V\cap V^\perp=0$.
Take $V$ to be a minimal $K$-submodule of $S\p$ generated by $S\a$.
Since $S\a$ is stable under complex conjugation, so is $V$,
because $V$ is spanned by the vectors of the form $k\cdot f$, $k\in K$,
$f\in S\a$, and since $\sigma$ commutes with complex conjugation on $\g$,
$K$ is stable under complex conjugation (we are in the adjoint group!).

By definition of Reynolds operator, $\Im\theta=V\cap (S\p)^K$.
By the construction of $B$ on $S\p$, it is clear that
$\Ker\psi=(S\a)^\perp\cap (S\p)^K$ and by $K$-invariance of $B$ we have
$\Ker\psi=V^\perp\cap (S\p)^K$. Since $V\cap V^\perp=0$, we finally obtain
$\Ker\psi\cap\Im\theta=0$.
\qed\enddemo

So to prove the surjectivity of $\psi$ it remains to prove that $\theta$
is injective. This idea (as well as many others in this paper)
belongs to Joseph, the only exception is that Joseph didn't
use Reynolds operator directly. Instead of that he defined
the ``Letzer map'' purely in infinitesimal terms, but it is easy
to see that in the adjoint situation Joseph's ``Letzer map''
actually coincides with good old Reynolds operator!

Now it is time to slightly change our notations and
to introduce new ones. First of all, since $S\p$ is a $K$-submodule of $S\g$,
injectivity of $\theta$ will follow from the injectivity of the restriction
on $(S\a)^{W_0}$ of the Reynolds operator $S\g\to(S\g)^K$.
Since the definition of the Reynolds operator depends only on the $K$-module
structure of $S\g$, we can substitute for $S\g$ the universal enveloping algebra
$U\g$, since $U\g$ has the same $K$-module structure by the
Poincare-Birkhoff-Witt theorem
(this $K$-module isomorphism could be written down explicitly as a
symmetrization map $S\g\to U\g$).
Moreover, we have an isomorphism $S\a\isom U\a$.
Notice that  the obvious embedding
$U\a\subset U\g$ corresponds to an embedding $S\a\subset S\g$ via the
symmetrization map.
So from now on $\theta$ will denote the map $(U\a)^{W_0}\to (U\g)^K$
induced by the Reynolds operator. We need to show that it is injective.

Fix a torus $S\subset G$ such that $\a=\Lie S$.
Then $S$ is a maximal $\sigma$-anisotropic torus in the terminology of
\cite{Vu}.
Fix a Borel subgroup $B$ and maximal torus $T\subset B$
such that $KB$ is dense in $G$
(so $K$ is spherical in $G$ and
if $\b$ is a Lie algebra of $B$ then $\k+\b=\g$), $T$ contains $S$
(such $B$, $T$ exist by \cite{Vu}).
Then automatically $\sigma(T)=T$.
Let $\h=\Lie T$
be a Cartan subalgebra of $\g$, let $\h_0=\h\cap\k$, so $\h=\h_0\oplus\a$.
Let ${\frak X}^*(T)$
(resp.~${\frak X}_*(T)$) be a character group
(resp.~the set of one-parameter subgroups of~$T$). We use an additive notation
for multiplication in ${\frak X}^*(T)$ and will sometimes identify a character
with its differential.
Let $P\subset{\frak X}^*(T)$ be the weight lattice,
let $P_+\subset P$ be the dominant weights.
Fix the sublattice $Q\subset P$ such that $\lambda\in Q$ if and only if
$$\sigma(\lambda)=-\lambda,\qquad \lambda|_S\in2{\frak X}^*(S).$$
It is clear that for any $\lambda\in Q$ we have $\lambda|_{\h_0}=0$,
so we may view such $\lambda$ as an element $\tilde\lambda$ of $\a^*$.
Let $w\in W_0$, $\tilde\mu=w\tilde\lambda$. Then we can extend $\tilde\mu$
to an element $\mu$ of ${\frak X}^*(T)\otimes{\Bbb Q}$. But it is clear that
actually $\mu\in Q$. Moreover by the structure theory of a baby Weyl group
(see eg.~\cite{Wa, \S1.1}), there exists $w\in W$ such that $\mu=w\lambda$.
Another useful tip is that if two elements of $\a^*$ are conjugate with
respect to $W$, then they are also conjugate by an element of $W_0$
(this follows from the fact that any $W_0$-fundamental chamber in $\a^*$
entirely lies in some $W$-fundamental chamber, see \cite{loc.~cit.}).
Set $Q_+=Q\cap P_+$.

\proclaim{Lemma 1.2}
In the identifications above, $Q=W_0Q_+$.
\endproclaim

\demo{Proof}
This follows from the results of T.~Vust (see \cite{Vu}).
Namely, passing to dual lattices,
it is sufficient to prove that if $R\in{\frak X}_*(S)$
is any one-parameter subgroup then there exists $w\in W_0$ such that
$R^w=wRw^{-1}$ is dominant, that is, for any $b\in B$ the limit
$\lim\limits_{t\to0}R^w(t)bR^w(t)^{-1}$ exists.
Consider a subgroup
$$P(R)=\{g\in G\,|\,\lim\limits_{t\to0}R(t)gR(t)^{-1}\ \text{exists}\}.$$
Then $P(R)$ is parabolic, moreover,
$P(R)$ and $\sigma P(R)$ are opposite parabolics, so that
$P(R)\cap\sigma P(R)=Z_G(R)$, the centralizer of $R$ in $G$.
Such parabolics are called $\sigma$-anisotropic.
Since $Z_G(R)$ is a $\sigma$-invariant reductive subgroup and $S\subset Z_G(R)$,
there exists minimal $\sigma$-anisotropic parabolic subgroup
$P'$ of $Z_G(R)$ such that $P'\cap\sigma P'=Z_{Z_G(R)}(S)$.
Then $P''=P'R_u(P(R))$ is a minimal $\sigma$-anisotropic
parabolic subgroup in $G$ and $P''\cap\sigma P''=Z_G(S)$.
Since $Z_G(R)\subset P(R)$, we have $P''\subset P(R)$.
Now, it is known that $K^0SR_u(P'')$ is open in~$G$, where $K^0$ is an identity
component of $K$
(algebraic analogue of Iwasawa decomposition), so if $B'$ is a Borel subgroup of
$Z_G(S)$ then
$B''=B'R_u(P'')$ is a Borel subgroup in $G$ that contains $S$,
is opposite to~$K$, and is contained in $P(R)$.
It remains to notice that all pairs $(B, T)$ such that
$B$ is a Borel subgroup opposite to~$K$, $T\subset B$ is a maximal
torus that contains a maximal $\sigma$-anisotropic torus,
are conjugate with respect to $K$ (see~\cite{loc. cit.}).
In particular, there exists $\tilde w\in N_K(S)$ such that
$B=\tilde wB''$. Since $B''\subset P(R)$, $B\subset P(\tilde wR\tilde w^{-1})$ and
hence
$R^w=\tilde wR\tilde w^{-1}$ is dominant.
\qed\enddemo

Returning to proof of Theorem,
for any spherical subgroup $H$ of $G$ and any irreducible representation
$V_\lambda$ of $G$ with highest weight $\lambda$, it is
known that $\dim V_\lambda^H\le 1$
and the set
$$\Gamma_H=\{\lambda\in P_+\,|\,\dim V_\lambda^H=1\}$$
is a submonoid in $P_+$ (see eg.~\cite{Pa2}).
In our case by \cite{Vu}, it is known that $\Gamma_K=Q_+$.
Irreducible representations $V_\lambda$ such that $\lambda\in Q_+$
are known as representations of class 1.
Fix $\lambda\in Q_+$,
for any $v\in V_\lambda$ and $\mu\in P$, let $(v)_\mu$ denote
a $\mu$-weight component of $v$ and
$\Supp v=\{\mu\in P\,|\,(v)_\mu\ne0\}$.
Let $v_\lambda$ be a highest weight vector
in~$V_\lambda$.
Fix
a $K$-invariant vector
$v_K\ne0$
in~$V_\lambda$.

\proclaim{Lemma 1.3}\par
{\rm(}i{\rm)} $\Supp v_K\subset Q${\rm;}\par
{\rm(}ii{\rm)} $\lambda\in\Supp v_K${\rm;}\par
{\rm(}iii{\rm)} $\Supp v_K\cap W\lambda=W_0\lambda$.
\endproclaim

\demo{Proof}
Fix $\mu\in\Supp v_K$.
Since $\h_0v_K=0$ we have $\mu(\h_0)=0$, and so $\sigma\mu=-\mu$.
Moreover, if $\tau\in S$ has order $2$, then $\tau\in K$,
so $\tau v_K=v_K$ and $\mu\in2{\frak X}^*(S)$. This proves (i).
To prove (ii) let us notice that if $B^0$ is a Borel subgroup, opposite
to $B$ (so that $B\cap B^0=T$), then $KB^0$ is dense in $G$
(this follows from the description of B as in proof of Lemma 1.2).
Then $\g=\b^0+\k$, so $U\g=U\b^0U\k$.
Then
$$V_\lambda=(U\g)v_K=(U\b^0)v_K.$$
But for any $x\in U\b^0$, $\Supp xv_K<\Supp v_K$
(with respect to a usual order relation on $P$),
so $\lambda\in \Supp v_K$.
The last assertion follows from (i), (ii), $N_K(S)$-invariance
of $v_K$ and the discussion before Lemma 1.2.
\qed\enddemo

Now we can finish the proof of the theorem in the case
of one summand. Suppose that $f\in\Ker\theta$, $f\ne0$.
To obtain a contradiction, it is sufficient to prove that
$f$ (viewed as an element of $S\a$) vanishes at $\tilde\lambda$
for all $\lambda\in Q$. Since $f$ is $W_0$-invariant
we need to prove that
$f(\tilde\lambda)=0$ for all $\lambda\in Q_+$
(by Lemma 1.2). Arguing by induction on the standard order relation
in $P_+$, we may assume that
$f(\tilde\mu)=0$
for all $\mu\in Q_+$, $\mu<\lambda$. Suppose that
$f(\tilde\lambda)\ne0$.
Let $\pi$ denote the $K$-invariant projection $V\lambda\to V_\lambda^K$.
Let $F$ denote the operator on $V_\lambda$ of left multiplication by $f$
(viewed as an element of $U\g$).
Then, clearly, $F(v_\lambda)=f(\tilde\mu)v_\lambda$.
Since $\theta f=0$ and since Reynolds operator commutes with
$K$-module homomorphisms, we have $\pi F\pi=0$
(because $\pi F\pi$ is clearly $K$-invariant, so Reynolds operator
maps it onto itself, on the other hand, Reynolds operator
commutes with left (and right) multiplication on ($K$-invariant!) element $\pi$).
Hence $\pi F v_K=0$. By Lemma 1.3 and by induction we have
$$Fv_K=\rho f(\tilde\lambda) \sum_{w\in W_0}\tilde w v_{\lambda},$$
where $\tilde w$ is some representative of $w$ in $N_K(S)$ and $\rho\ne0$.
Since $\pi$ is $K$-invariant, we get
$\pi(v_\lambda)=0$. Therefore
$(U\k)v_\lambda\ne V_\lambda$. But
$$V_\lambda=(U\g)v_\lambda=(U(\k+\b))v_\lambda=(U\k)(U\b)v_\lambda=(U\k)v_\lambda.$$
Contradiction.

\head\S2. Surjectivity of
$\psi:\,\C[\p\times \p]^K\to\C[\a\times\a]^{W_0}$\endhead
First, we can identify
$C[\p\times \p]$ with
$S(\p\times \p)$, and
$C[\a\times \a]$ with
$S(\a\times \a)$ by means of Killing form.
Then we can identify
$S(\a\times \a)$ with
$U(\a\times \a)$.
Let $\theta$ denote a restriction on
$U(\a\times \a)^{W_0}$ of a Reynolds operator
$U(\g\times\g)\to U(\g\times\g)^K$
(with respect to diagonal action).

\proclaim{Lemma 2.1}
If $\theta$ is injective then $\psi$ is surjective.
\endproclaim

\demo{Proof}
The proof is the same as of Lemma 1.1 and the discussion after it.
\qed\enddemo

Now suppose that $f\in\Ker\theta$. By Lemma 1.2
we know that $Q_+\times Q_+$ is Zarisky dense in $\C\otimes(Q_+\times Q_+)$.
Hence to show that $f=0$
it is sufficient to show that $f$ (viewed as an element of
$S(\a\times\a)$) vanishes at all $(\lambda,\mu)\in Q_+\times Q_+$.
We will proceed by a usual induction on the standard order relation.
So suppose that $f$ vanishes at all $(\lambda',\mu')<(\lambda,\mu)$.
Consider the $G\times G$ module
$V_\lambda\otimes V_\mu$.
Let $\pi$ denote a $K$-invariant projection onto
$(V_\lambda\otimes V_\mu)^K$ (with respect to a diagonal action).
Let $F$ denote the operator of left multiplication by $f$
(viewed as an element of $U(\g\times\g)$).
Since $\theta f=0$, we get $\pi F\pi=0$ as in the proof
for a case of one summand. Let $v_K^\lambda$ (resp.~$v_K^\mu$)
denote a non-zero $K$-invariant vector in $V_\lambda$ (resp.~$V_\mu$).
Then $\pi F(v_K^\lambda\otimes v_K^\mu)=0$.
Let $v_\lambda$ (resp.~$v_\mu$)
denote the highest weight vector of $V_\lambda$ (resp. $V_\mu$).
By Lemma 1.3
$$v_K^\lambda\otimes v_K^\mu=
\sum_{{w'\in W_0/(W_0)_\lambda\atop
w''\in W_0/(W_0)_\mu}}(\tilde w'v_\lambda)\otimes (\tilde w''v_\mu)+U,$$
where $\Supp U$ lies in $Q\times Q$ and for any
$(\lambda',\mu')\in\Supp U$ by induction we have $f(\lambda',\mu')=0$.
So we obtain that
$$\sum_{{w'\in W_0/(W_0)_\lambda\atop
w''\in W_0/(W_0)_\mu}}f(w'\lambda,w''\mu)(\tilde w'v_\lambda)\otimes
(\tilde w''v_\mu)$$
lies in the kernel of $\pi$.
Now for any
$w'\in W_0/(W_0)_\lambda$, $w''\in W_0/(W_0)_\mu$
we can find an element $\left(w_1(w',w''), w_2(w',w'')\right)\in W_0\times W_0$
such that $w_1=yw'$, $w_2=yw''$ for some $y\in W_0$
and such that $w_1\lambda+w_2\mu\in Q_+$.
Denote by $w_1^i$, $w_2^i$, $i=1,\ldots,r$ the complete set of
such pairs obtained for all
$w'\in W_0/(W_0)_\lambda$, $w''\in W_0/(W_0)_\mu$.
Now by $W_0$-invariance of $f$ and by $K$-invariance of $\pi$ we
proved that an element
$$\sum_{i=1,\ldots r}
a_if(w_1^i\lambda,w_2^i\mu)
(\tilde w_1^iv_\lambda)\otimes
(\tilde w_2^iv_\mu)$$
lies in the kernel of $\pi$ for some non-zero integers $a_r$.
It is sufficient to prove that all
$f(w_1^i\lambda,w_2^i\mu)$ are equal to zero.
So we reduced our proof to the following claim:

\proclaim{Claim}
Suppose that an element
$$A=\sum_{i=1,\ldots r}
b_i
(\tilde w_1^iv_\lambda)\otimes(\tilde w_2^iv_\mu)$$
lies in the kernel of $\pi$. Then all $b_i$ are equal to zero.
\endproclaim

Suppose that some $b_i\ne0$.
By Kostant's refinement of PRV conjecture (see \cite{J} or \cite{Ku1})
for any $i$ the module $(U\g)\left[\tilde w_1^iv_\lambda)\otimes(\tilde
w_2^iv_\mu)\right]$
containes a unique $U\g$-submodule of highest weight
$w_1^i\lambda+w_2^i\mu$.
Denote it heighest weight vector by $v^i$.
By Kumar's refinement of PRV conjecture (see \cite{J} or \cite{Ku2}),
these $U\g$ submodules form a direct sum
$S=\mathop{\oplus}\limits_{i=1}^rS_i$.
Set
$$S'=\mathop{\oplus}_{{\scriptstyle 1\le i\le r}\atop{\scriptstyle b_i\ne
0}}S_i.$$
By Density theorem it follows that the module $(U\g) A$
contains $S'$.
Let $\pi_0$ denote a $G$-invariant projection $(U\g) A\to S$.
Take $i$ such that the weight
$w_1^i\lambda+w_2^i\mu$ is minimal.
Let $\pi_1$ denote a $G$-invariant projection $S'\to S_i$.
Then it is clear from the above and by checking the weights that
$$\pi_1\pi_0(A)=cv^i,$$
where $c\ne0$. Denote by $\pi_2$ a $K$-invariant
projection $S_i\to S_i^K$.
Then arguing as in proof for one summand we see
that $\pi_2\pi_1\pi_0(A)\ne0$. But then $\pi(A)\ne0$ as well,
because $\pi_2\pi_1\pi_0$ is $K$-invariant and hence factors through
$\pi|_{(U\g) A}$.
This completes the proof.

\enddocument